# On Affine and Conjugate Nonparametric Regression

BY RAJESHWARI MAJUMDAR

*University of Connecticut*

Suppose the nonparametric regression function of a response variable $Y$ on covariates $X$ and $Z$ is an affine function of $X$ such that the slope $\beta$ and the intercept $\alpha$ are real valued measurable functions on the range of the completely arbitrary random element $Z$. Assume that $X$ has a finite moment of order greater than or equal to 2, $Y$ has a finite moment of conjugate order, and $\alpha(Z)$ and $\alpha(Z)X$ have finite first moments. Then, the nonparametric regression function equals the least squares linear regression function of $Y$ on $X$ with all the moments that appear in the expression of the linear regression function calculated conditional on $Z$. Consequently, conditional mean independence implies zero conditional covariance and a degenerate version of the aforesaid affine form for the nonparametric regression function, whereas the aforesaid affine form and zero conditional covariance imply conditional mean independence. Further, it turns out that the nonparametric regression function has the aforesaid affine form if $X$ is Bernoulli, and since 1 is the conjugate exponent of $\infty$, the least squares linear regression formula for the nonparametric regression function holds when $Y$ has only a finite first moment and $Z$ is completely arbitrary.

**1. Introduction.** The nonparametric regression function of a response variable $Y$ on covariates $X$ and $Z$, that is, the conditional expectation of $Y$ given $X$ and $Z$, where $X$ is a Bernoulli random variable, is frequently of interest in the social sciences. A couple of specific instances are provided by the study of microeconometric treatment effects [Imbens and Wooldridge (2007)] and by the quantification of the spatial competition model [Smithies (1941)] in the context of analyzing the statics and dynamics of party ideologies in a two-party electoral system [Downs (1957)]; we will discuss these examples in greater detail later, but for now, we concentrate on obtaining a tractable expression for the nonparametric regression function.

Given random variables $X$ and $Y$, and a random element $Z$, let

$$g(X,Z) = \frac{\mathrm{Cov}^Z(X,Y)}{\mathrm{Var}^Z(X)}\big[\mathrm{Var}^Z(X) > 0\big]\big(X - \mathrm{E}^Z(X)\big), \qquad (1)$$

where the conditional expectation of a random variable $\mathcal{G}$ given a random element $\mathfrak{N}$ is denoted by $\mathrm{E}^{\mathfrak{N}}(\mathcal{G})$, and an event $A$ is identified with its indicator function. We show in

---

this note that the formula

$$\mathrm{E}^{X,Z}(Y) = \mathrm{E}^{Z}(Y) + g(X, Z) \qquad (2)$$

holds if $X$ is Bernoulli, under the minimal assumption that $Y$ has a finite first moment, with $Z$ completely arbitrary.

While it may be tempting to try to verify (2) for Bernoulli $X$ from the definition of conditional expectation, we doubt if it can be done. Clearly, $\mathrm{E}^{Z}(Y) + g(X, Z)$, the right-hand side of (2), is a measurable function of $(X, Z)$; if it were to equal $\mathrm{E}^{X,Z}(Y)$, then it must be integrable, since $\mathrm{E}^{X,Z}(Y)$ itself is integrable [Dudley (1989, Theorem 10.1.1)]. Since $\mathrm{E}^{Z}(Y)$ is integrable as well, the first order of business in establishing (2) would be to verify that $g(X, Z)$ is integrable. As discussed in Remark 7, the verification of the integrability of $g(X, Z)$ for Bernoulli $X$ is, for all practical intents and purposes, impossible unless one is willing to assume that $Y$ is bounded, or assume that $Y$ has a finite second moment and the range of $Z$ is a Polish space.

We obtain (2) for Bernoulli $X$ from Theorem 1, which results from our investigation of the hierarchy of the interdependence among the notions of conditional independence, conditional mean independence, and zero conditional covariance. It is well-known that conditional independence implies both conditional mean independence and zero conditional covariance. We present counterexamples to show that the direction of neither of these hierarchies can be reversed. It is also well-known that conditional mean independence implies zero conditional covariance, but the reverse implication does not hold in general (Remark 5), though it does for Bernoulli $X$ or jointly Normal $(X, Y, Z)$.

That the reverse implication holds if $X$ is Bernoulli is commonly proved using the assumption that the propensity score [Rosenbaum and Rubin (1983)] $\mathrm{E}^{Z}(X)$ is non-degenerate, that is, in the open interval $(0, 1)$ [Imbens and Wooldridge (2007)], so that division by it is legitimate, and the fact that $X = [X = 1]$. On the other hand, that the reverse implication holds if $(X, Y, Z)$ is multivariate Normal is commonly proved by using the fact that zero conditional covariance in this context implies conditional independence, and hence conditional mean independence. These context specific arguments obscure the common structural thread that connects the two cases; to identify that connecting thread, we offer below a different take on the reverse implication in the multivariate Normal case.

If $(X, Y, Z)$ is multivariate Normal, such that $Z = (Z_1, \cdots, Z_{n-2}) \in \Re^{n-2}$ and $(X, Y) \in \Re^{2}$, then the formula in (2) holds. We believe that to be the core reason why zero conditional covariance implies conditional mean independence in the multivariate Normal case. That (2) holds in the multivariate Normal case can be verified by directly computing the various conditional expectations that appear in (2), but the computations are tedious if the covariance matrix of $(X, Y, Z)$ is invertible, and a nightmare otherwise, as they then involve calculation of generalized inverses of partitioned matrices; the orthogonality argument outlined below provides an elegant alternative.

By Theorem 1.2.11 of Muirhead (1982), the conditional distributions of $Y$ given $(X, Z) \in \Re^{n-1}$, of $Y$ given $Z \in \Re^{n-2}$, and of $X$ given $Z \in \Re^{n-2}$ are all Normal on the line, the conditional distribution of $(Y, X)$ given $Z \in \Re^{n-2}$ is Normal on the plane, the mean of each of these distributions is an affine function of the conditioning element in question, and the covariance structure of none of these distributions depends on the conditioning element in question; that is,

$$\mathrm{E}^Z(Y) = d + \sum_{j=1}^{n-2} e_i Z_i \text{ for } d, e_1, \cdots, e_{n-2} \in \Re$$

$$\mathrm{E}^Z(X) = k + \sum_{j=1}^{n-2} h_i Z_i \text{ for } k, h_1, \cdots, h_{n-2} \in \Re$$

$$\mathrm{E}^{X,Z}(Y) = a + bX + \sum_{j=1}^{n-2} c_i Z_i \text{ for } a, b, c_1, \cdots, c_{n-2} \in \Re \quad (3)$$

$$\mathrm{Cov}^Z(X, Y) = \mathrm{E}\big((X - \mathrm{E}^Z(X))(Y - \mathrm{E}^Z(Y))\big)$$

$$\mathrm{Var}^Z(X) = \mathrm{E}\big((X - \mathrm{E}^Z(X))^2\big).$$

Since $\mathrm{E}^Z(Y) = \mathrm{E}^Z\big(\mathrm{E}^{X,Z}(Y)\big)$, we obtain

$$\mathrm{E}^{X,Z}(Y) - \mathrm{E}^Z(Y) = b\big(X - \mathrm{E}^Z(X)\big), \quad (4)$$

and (2) will follow once we show that

$$b = \frac{\mathrm{Cov}^Z(X, Y)}{\mathrm{Var}^Z(X)} \big[\mathrm{Var}^Z(X) > 0\big]. \quad (5)$$

In the Hilbert space $H$ of square integrable functions on the probability space underlying $(X, Y, Z)$, let $V^\perp$ denote the orthogonal complement of a closed subspace $V$ and $\Pi_V$ the orthogonal projection operator on to $V$. Let $\mathfrak{L}(X, Z)$ (respectively, $\mathfrak{L}(Z)$) denote the finite-dimensional subspace of $H$ spanned by $\{J, X, Z\}$ (respectively, $\{J, Z\}$), where $J \equiv 1$. It is well-understood that $\mathfrak{L}(X, Z)$ (respectively, $\mathfrak{L}(Z)$) is contained in a closed subspace $\mathfrak{M}(X, Z)$ (respectively, $\mathfrak{M}(Z)$) of $H$ so that

$$\mathrm{E}^{X,Z}(Y) = \Pi_{\mathfrak{M}(X,Z)} Y, \ \mathrm{E}^Z(Y) = \Pi_{\mathfrak{M}(Z)} Y, \text{ and } \mathrm{E}^Z(X) = \Pi_{\mathfrak{M}(Z)} X;$$

see Remark 1 of Majumdar (2017) for an outline of the argument validating the above. Note that, by (3), $\mathrm{E}^{X,Z}(Y) \in \mathfrak{L}(X, Z)$, $\mathrm{E}^Z(Y) \in \mathfrak{L}(Z)$, and $\mathrm{E}^Z(X) \in \mathfrak{L}(Z)$; thus, since $\mathfrak{L}(X, Z) \subset \mathfrak{M}(X, Z)$,

$$\mathrm{E}^{X,Z}(Y) = \Pi_{\mathfrak{L}(X,Z)} \mathrm{E}^{X,Z}(Y) = \Pi_{\mathfrak{L}(X,Z)} \Pi_{\mathfrak{M}(X,Z)} Y = \Pi_{\mathfrak{L}(X,Z)} Y,$$

and, since $\mathfrak{L}(Z) \subset \mathfrak{M}(Z)$,

$$\mathrm{E}^Z(Y) = \Pi_{\mathfrak{L}(Z)} \mathrm{E}^Z(Y) = \Pi_{\mathfrak{L}(Z)} \Pi_{\mathfrak{M}(Z)} Y = \Pi_{\mathfrak{L}(Z)} Y$$

$$\mathrm{E}^Z(X) = \Pi_{\mathfrak{L}(Z)} \mathrm{E}^Z(X) = \Pi_{\mathfrak{L}(Z)} \Pi_{\mathfrak{M}(Z)} X = \Pi_{\mathfrak{L}(Z)} X.$$

Further, since the conditional variance and the conditional covariance in (3) do not depend on the conditioning element,

$$\mathrm{Cov}^Z(X, Y) = \langle \Pi_{(\mathfrak{L}(Z))^\perp} Y, \Pi_{(\mathfrak{L}(Z))^\perp} X \rangle \text{ and } \mathrm{Var}^Z(X) = \left\| \Pi_{(\mathfrak{L}(Z))^\perp} X \right\|^2.$$

Clearly, if $X \in \mathfrak{L}(Z)$, then $\mathrm{Var}^Z(X) = 0$, rendering RHS(5) = 0, and we can take $b = 0$ in (3) through a reparametrization. Otherwise, since

$$\mathfrak{L}(X, Z) = \{V + cX : V \in \mathfrak{L}(Z), c \in \Re\},$$

and RHS(5) $= \left\| \Pi_{(\mathfrak{L}(Z))^\perp} X \right\|^{-2} \langle \Pi_{(\mathfrak{L}(Z))^\perp} Y, \Pi_{(\mathfrak{L}(Z))^\perp} X \rangle$, (4) follows from Lemma 1 (the "partialling out" lemma) below.

**Lemma 1** For a proper and closed subspace $V$ of a Hilbert space $H$ and $x \in H \setminus V$, let $V_x$ denote the closure of the subspace $\{v + cx : v \in V, c \in \Re\}$. Then, for any $y \in H$,

$$\Pi_{V_x} y = \Pi_V y + \left\| \Pi_{V^\perp} x \right\|^{-2} \langle \Pi_{V^\perp} y, \Pi_{V^\perp} x \rangle \Pi_{V^\perp} x.$$

<u>Proof of Lemma 1.</u> By the definition of $V_x$, there exists a sequence $\{v_n : n \geq 1\} \subset V$ and a sequence $\{c_n : n \geq 1\} \subset \Re$ such that

$$\Pi_{V_x} y = \lim_{n \to \infty} (v_n + c_n x). \tag{6}$$

Since $V \subset V_x$ and the orthogonal projection operator is continuous,

$$\Pi_V y = \lim_{n \to \infty} (v_n + c_n \Pi_V x),$$

implying, by (6) and the definition of the orthogonal projection operator,

$$\Pi_{V_x} y - \Pi_V y = \lim_{n \to \infty} c_n \Pi_{V^\perp} x. \tag{7}$$

By continuity and homogeneity of inner product and (7), since $x \notin V$,

$$\lim_{n \to \infty} c_n = \left\| \Pi_{V^\perp} x \right\|^{-2} \langle \Pi_{V_x} y - \Pi_V y, \Pi_{V^\perp} x \rangle. \tag{8}$$

Since $x \in V_x$ and $\Pi_V x \in V \subset V_x$, we obtain

$$\Pi_{V^\perp} x \in V_x;$$

since $\Pi_{V_x} y - \Pi_V y = \Pi_{V^\perp} y - \Pi_{(V_x)^\perp} y$, the lemma follows from (8) and (7). □

Now note that, if $X$ is Bernoulli then

$$\mathrm{E}^{X,Z}(Y) = e(X, Z) = e(0, Z) + (e(1, Z) - e(0, Z))X.$$

Thus, the representation

$$\mathrm{E}^{X,Z}(Y) = \alpha(Z) + \beta(Z)X \tag{9}$$

holds for measurable functions $\alpha$ and $\beta$ on the range of $Z$ in the Bernoulli $X$ case (as well as in the multivariate Normal case, where the slope does not depend on $Z$). It is quite conceivable that the reverse implication holds in the Bernoulli $X$ case because the affine structure in (9) relates the slope to the conditional covariance the way it does in the multivariate Normal case. That leads to the conjecture that (9), under suitable (but fairly general) assumptions (that are satisfied when $X$ is Bernoulli), implies (2). Theorem 1 settles this conjecture by showing that if $X$ has a finite moment of order greater than or equal to 2, $Y$ has a finite moment of conjugate order, and (9) holds such that $\alpha(Z)$ and $\alpha(Z)X$ have finite first moments, then (2) holds.

We now return to elaborating the two examples we mentioned earlier on the importance of a Bernoulli covariate in the social sciences. Imbens and Wooldridge (2007) cite the problem of measuring (by $Y$) the effectiveness of a particular training protocol in a labor market program, where some individuals receive the training and others do not (with $X$ tracking the status of the individuals vis-a-vis the receipt of the training), and each subject has a vector of characteristics (represented by $Z$). Typically, one is interested in the estimation of either the population average treatment effect (PATE), defined as

$$\mathrm{E}(Y[X=1] - Y[X=0]),$$

or the population average treatment effect for the treated (PATT), defined as

$$\mathrm{E}^{[X=1]}((Y[X=1] - Y[X=0])) = (2[X=1] - 1)\mathrm{E}^{[X=1]}(Y).$$

However, if one is interested in the estimation of the treatment effect controlling for the covariate that represents the individual characteristics, the expression of $\mathrm{E}^{X,Z}(Y)$ obtained in (2) allows one to recognize the treatment effect as the difference

$$g(1, Z) - g(0, Z) = \frac{\mathrm{Cov}^Z(X, Y)}{\mathrm{Var}^Z(X)} [\mathrm{Var}^Z(X) > 0];$$

if the propensity score is assumed to be non-degenerate, then $[\mathrm{Var}^Z(X) > 0] = 1$, reducing the treatment effect to the ratio $\mathrm{Cov}^Z(X, Y)/\mathrm{Var}^Z(X)$.

In a two-party electoral system the strategies employed by the parties are, to a large extent, determined by the distribution of voters on the ideological orientation scale [Downs (1957, page 117)]. Conceptually, if the distribution is unimodal and symmetric, the parties are incentivized to converge towards the center, whereas if the distribution is bimodal with the two modes situated at the two ends of the ideological spectrum, the parties have very little motivation to modify their platforms to catch the voters in the middle. If the distribution is unimodal and skewed to the right (respectively, left), then the party with a platform tilted towards the left (respectively, right) end of the spectrum has an edge to start with (which is difficult to nullify over the course of an election cycle). Since each party has the ultimate objective of winning the election, the outcome of the election, denoted by the Bernoulli variable $Y$, where

$$Y = \begin{cases} 1 & \text{if party } A \text{ wins the election} \\ 0 & \text{otherwise,} \end{cases}$$

is a response variable of interest. The distribution of voters on the ideological orientation scale, denoted by $Z$, is an important functional covariate for $Y$. Another important covariate is the quality and intensity of the campaign by the parties, measured by how successfully the parties adapt their platforms to the ideological orientation scale in their ambient environment and communicate their evolving positions to the voters; this covariate can be dichotomized and denoted by $X$, where

$$X = \begin{cases} 1 & \text{if party } A \text{ runs a more effective campaign} \\ 0 & \text{otherwise.} \end{cases}$$

Thus, given the distribution of voters on the ideological orientation scale, the problem of estimating how the probability of electoral success for party $A$ will change as its campaign becomes more effective is reduced, by (2), to estimating the ratio of the conditional covariance to the conditional variance (of course, under the assumption that the conditional probability of party $A$ running a more effective campaign is never 0 or 1). It is possible that for $Z$ in the neighborhood of symmetric and unimodal distributions, the probability of electoral success will not depend on the campaign quality, but for other values of $Z$, it will. This regression formulation of the problem creates a pathway to formally testing that possibility.

The following notational conventions are used throughout the note. Equality (or inequality) involving measurable functions defined on a probability space indicates the relation holds almost surely. The universal null set is denoted by $\mathfrak{K}$ and the Borel $\sigma$-algebra of a metric space $\mathfrak{Z}$ by $\mathcal{B}(\mathfrak{Z})$. The notation $\otimes$ denotes the product of $\sigma$-algebras. The Normal distribution with mean $\mu$ and variance $\sigma^2$ is denoted by $\mathcal{N}(\mu, \sigma^2)$ and the density of the $\mathcal{N}(0,1)$ distribution by $\psi$. For $\nu \in (-1, 1)$, $\Phi_\nu$ denotes the bivariate normal measure on $\mathcal{B}(\Re^2)$ with means 0, variances 1, and covariance $\nu$, and $\phi_\nu$ the density of $\Phi_\nu$. We also repeatedly use the *averaging property*, the *chain rule*, and the *pull-out property* of conditional expectation [Kallenberg (2002, page 105)], and fastidiously verify the underlying (minimal) integrability assumptions before using these very well-known results.

**2. Results.** Throughout the rest of the note, we will work with the following general structure, which, as necessary, will be appropriately restricted.

Let $(\Omega, \mathfrak{F}, P)$ be a probability space and E the induced expectation. Let $X, Y, Z$ be random elements defined on $(\Omega, \mathfrak{F}, P)$, where the range of $X$ (respectively, $Y$) is a Polish space $\mathfrak{X}$ (respectively, $\mathfrak{Y}$) endowed with its Borel $\sigma$-algebra $\mathcal{B}(\mathfrak{X})$ (respectively, $\mathcal{B}(\mathfrak{Y})$) and the range of $Z$ is an arbitrary measurable space $(\mathfrak{S}, \mathfrak{T})$. Consistent with the notation for conditional expectation given a random element introduced in Section 1, let $\mathrm{E}^{\mathfrak{C}}$ denote the conditional expectation given a sub $\sigma$-algebra $\mathfrak{C}$ of $\mathfrak{F}$.

**Definition 1.** The *conditional distribution* of $X$ given $\mathfrak{C}$, denoted by $P_{X|\mathfrak{C}}(\cdot, \cdot)$, is a function on $\mathcal{B}(\mathfrak{X}) \times \Omega$ taking values in $[0, 1]$ such that

for $P$-almost all $\omega$, $P_{X|\mathfrak{C}}(\,\cdot\,,\omega)$ is a probability measure on $\mathcal{B}(\mathfrak{X})$, and for each $A \in \mathcal{B}(\mathfrak{X})$, $P_{X|\mathfrak{C}}(A,\,\cdot\,) = \mathrm{E}^{\mathfrak{C}}([X \in A])(\,\cdot\,)$;

see Dudley (1989, Sections 10.1 and 10.2). $P_{Y|\mathfrak{C}}$ and $P_{(X,Y)|\mathfrak{C}}$ are similarly defined. Note that $P_{X|\mathfrak{C}}$, $P_{Y|\mathfrak{C}}$, and $P_{(X,Y)|\mathfrak{C}}$ exist on $\mathcal{B}(\mathfrak{X}) \times \Omega$, $\mathcal{B}(\mathfrak{Y}) \times \Omega$, and $\mathcal{B}(\mathfrak{X} \times \mathfrak{Y}) \times \Omega$, respectively. Since $\mathfrak{X}$ and $\mathfrak{Y}$ are assumed to be Polish spaces, $\mathcal{B}(\mathfrak{X} \times \mathfrak{Y})$ equals $\mathcal{B}(\mathfrak{X}) \otimes \mathcal{B}(\mathfrak{Y})$, and the existence of the three conditional distributions is guaranteed by Theorem 10.2.2 of Dudley (1989).

**Definition 2.** The random elements $X$ and $Y$ are *conditionally independent* given $Z$, if, for every $A \in \mathcal{B}(\mathfrak{X})$ and $B \in \mathcal{B}(\mathfrak{Y})$,

$$P_{(X,Y)|Z}(A \times B,\,\cdot\,) = P_{X|Z}(A,\,\cdot\,) \times P_{Y|Z}(B,\,\cdot\,),$$

where $Z$ is the shorthand for $\sigma(Z)$, the $\sigma$-algebra generated by $Z$. Clearly, $X$ and $Y$ are conditionally independent given $Z$ means that the conditional distribution of $(X,Y)$ given $Z$ on $\mathcal{B}(\mathfrak{X}) \otimes \mathcal{B}(\mathfrak{Y})$ is the product measure of the conditional distribution of $X$ given $Z$ on $\mathcal{B}(\mathfrak{X})$ and the conditional distribution of $Y$ given $Z$ on $\mathcal{B}(\mathfrak{Y})$.

Next in line is the definition of conditional mean independence. In econometric analysis, Heckman et. al. (1998) has presented the notion of conditional mean independence with Bernoulli $X$ as an alternative to the notion of unconfoundedness [Rubin (1990)].

For an arbitrary probability space $(\Theta, \mathfrak{A}, Q)$ and $1 \leq p < \infty$, let $\mathcal{L}_p(\Theta, \mathfrak{A}, Q)$ denote the Banach space of real valued functions $f$ on $(\Theta, \mathfrak{A}, Q)$ such that $\mathrm{E}_Q(|f|^p) < \infty$, where $\mathrm{E}_Q$ is the expectation induced by $Q$; let $\mathcal{L}_\infty(\Theta, \mathfrak{A}, Q)$ denote the Banach space of essentially bounded functions on $(\Theta, \mathfrak{A}, Q)$ [Rudin (1986, Definition 3.7 and page 96)]. For $1 \leq p \leq \infty$, we write $\mathcal{L}_p$ for $\mathcal{L}_p(\Theta, \mathfrak{A}, Q)$ when there is no contextual ambiguity regarding the underlying probability space; also, $\|\cdot\|_p$ denotes the norm in $\mathcal{L}_p$.

**Definition 3.** When $\mathfrak{Y} = \mathfrak{R}$ and $Y \in \mathcal{L}_1$, $Y$ is *mean independent* of $X$ given $Z$ if

$$\mathrm{E}^{X,Z}(Y) = \mathrm{E}^Z(Y). \tag{10}$$

Similarly, when $\mathfrak{X} = \mathfrak{R}$ and $X \in \mathcal{L}_1$, $X$ is *mean independent* of $Y$ given $Z$ if

$$\mathrm{E}^{Y,Z}(X) = \mathrm{E}^Z(X). \tag{11}$$

For future use, we now present the Moment inequality of Lyapunov, a consequence of Corollary 1.4.17 of Fabian and Hannan (1985). Given a real valued measurable function $f$ on $(\Theta, \mathfrak{A}, Q)$, the function $u \mapsto \left[\mathrm{E}_Q(|f|^u)\right]^{1/u}$, defined on $[1, \infty]$, is non-decreasing. Note that $\left[\mathrm{E}_Q(|f|^u)\right]$ can be $\infty$ for all $u \geq 1$.

**Definition 4.** Let $1 \leq p, q \leq \infty$ be a pair of conjugate exponents [Rudin (1986, Definition 3.4)]. If $\mathfrak{X} = \mathfrak{Y} = \mathfrak{R}$, $X \in \mathcal{L}_p$, and $Y \in \mathcal{L}_q$, the *conditional covariance* of $X$

and $Y$ given $Z$ is

$$\text{Cov}^Z(X,Y) = \text{E}^Z\big((X - \text{E}^Z(X))(Y - \text{E}^Z(Y))\big) \\ = \text{E}^Z(XY) - \text{E}^Z(X)\text{E}^Z(Y). \quad (12)$$

By Holder's inequality [Rudin (1986, Theorem 3.8)],

$$XY \in \mathcal{L}_1. \quad (13)$$

Recall that the conditional distribution $P_{X|Z}$ exists on $\mathcal{B}(\Re)$ and if the conditional distribution exists, then the conditional expectation is the integral with respect to that distribution [Dudley (1989, Theorem 10.2.5)]. By the conditional Jensen's inequality [Dudley (1989, Theorem 10.2.7)] and the Moment inequality in $(\Re, \mathcal{B}(\Re), P_{X|Z})$

$$|\text{E}^Z(X)| \leq \text{E}^Z(|X|) \leq \big(\text{E}^Z(|X|^p)\big)^{1/p},$$

whence, by the averaging property, $\text{E}^Z(X) \in \mathcal{L}_p$. Similarly, $\text{E}^Z(Y) \in \mathcal{L}_q$, implying, by Holder's inequality, $(X - \text{E}^Z(X))(Y - \text{E}^Z(Y)) \in \mathcal{L}_1$; that makes $\text{Cov}^Z(X,Y)$ well-defined. Since each of the three random variables $Y\text{E}^Z(X)$, $X\text{E}^Z(Y)$, and $\text{E}^Z(Y)\text{E}^Z(X)$ is in $\mathcal{L}_1$ by Holder's inequality, the second equality in (12) follows from (13), the linearity of conditional expectation, and the pull-out property.

**Definition 5.** If $\mathfrak{X} = \Re$ and $X \in \mathcal{L}_p$ for $p \geq 2$, the *conditional variance* of $X$ given $Z$ is

$$\text{Var}^Z(X) = \text{E}^Z\big((X - \text{E}^Z(X))^2\big) = \text{E}^Z(X^2) - (\text{E}^Z(X))^2. \quad (14)$$

Since $p \geq 2$, the conjugate exponent $q \leq 2$. As such, by the Moment inequality, $X \in \mathcal{L}_q$ as well. Thus, by (12), $\text{Var}^Z(X)$ is well-defined and the second equality holds.

Proposition 1 presents a well-known characterization of conditional independence which is used to show that conditional independence implies conditional mean independence.

**Proposition 1.** The following statements are equivalent:

(1) $X$ and $Y$ are conditionally independent given $Z$.
(2) $P_{Y|(X,Z)} = P_{Y|Z}$.
(3) $P_{X|(Y,Z)} = P_{X|Z}$.

**Remark 1.** The equivalence of statements (1) and (2) in Proposition 1 is the assertion of Proposition 6.6 of Kallenberg (2002); interchange of $X$ and $Y$ establishes the equivalence of statements (1) and (3). The assertion of statement (3) in Proposition 1 for Bernoulli $X$ is what Rubin (1990) refers to as $X$ being unconfounded with $Y$ given $Z$. //

**Remark 2.** Suppose $X$ and $Y$ are conditionally independent given $Z$. Applying the Fubini Theorem in $(\Re^2, \mathcal{B}(\Re^2), P_{(X,Y)|Z})$ and (12), we obtain

$$\text{Cov}^Z(X,Y) = 0. \quad (15)$$

If $Y \in \mathcal{L}_1$, then the equivalence of statements (1) and (2) in Proposition 1 implies that (10) holds. Similarly, if $X \in \mathcal{L}_1$, then the equivalence of statements (1) and (3) in Proposition 1 implies that (11) holds. //

**Remark 3.** Unsurprisingly, (10) can hold without (11) holding. Recall from Example 2 of Majumdar (2017) that if we define random variables $X$ and $Y$ on the equi-probable discrete sample space $\Omega = \{-1, 0, 1\}$ as $X(\omega) = [\omega = 0]$ and $Y(\omega) = \omega$, then $\mathrm{E}^X(Y) = \mathrm{E}(Y)$ but $\mathrm{E}^Y(X) \neq \mathrm{E}(X)$. If we define $Z$ to be a degenerate random variable on this sample space so that $X$ and $Y$ are independent of $Z$, we have a trivial example of the notion of conditional mean independence being asymmetric. Interestingly, there is no non-degenerate random variable $Z$ on this sample space such that $X$ and $Y$ are independent of $Z$. In fact, as argued below, it is not possible to construct a non-degenerate $Z$ on this sample space such that (10) holds but (11) does not.

Suppose such a $Z$ exists. Since $X$ is a measurable function of $Y$, for every $Z$, $\mathrm{E}^{Y,Z}(X) = X$; for (11) to not hold, we must have

$$\mathrm{E}^Z(X) \neq X, \tag{16}$$

that is, $X$ cannot be a measurable function of $Z$. If $Z(0) = a$ differs from both $Z(1)$ and $Z(-1)$, then $X = [Z = a]$, and $X$ becomes a measurable function of $Z$; thus (16) implies that one of $Z(1)$ and $Z(-1)$ must equal $a$. Since $Z$ has to be non-degenerate, both $Z(1)$ and $Z(-1)$ cannot equal $a$. Assume, without loss of generality, that $Z(1) = a \neq Z(-1) = b$. Note that

$$\mathrm{E}^{X,Z}(Y) = 1 - X - 2[Z = b] = [X = 0] - 2[Z = b], \tag{17}$$

implied by the string of self explanatory equalities below:

$$\mathrm{E}(([X = 0] - 2[Z = b])A) = 0 = \mathrm{E}(YA) \quad \text{if } A = \mathfrak{K}, \Omega, \{0\}, \{-1, 1\}$$

$$\mathrm{E}(([X = 0] - 2[Z = b])A) = \frac{1}{3} = \mathrm{E}(YA) \quad \text{if } A = \{1\}, \{0, 1\},$$

$$\mathrm{E}(([X = 0] - 2[Z = b])A) = -\frac{1}{3} = \mathrm{E}(YA) \quad \text{if } A = \{-1\}, \{0, -1\}.$$

Now (17) implies, by the chain rule, $\mathrm{E}^Z(Y) = 1 - \mathrm{E}^Z(X) - 2[Z = b]$, which, by (17) and (16), violates the premise that (10) holds. //

Example 1 unveils $(X, Y, Z)$ where $Z$ is non-degenerate but (10) does not imply (11).

**Example 1.** For $\rho \in (-1, 1)$, let $(Y, Z)$ be the random vector in $\Re^2$ such that the conditional joint distribution of $(Y, Z)$ given $W$ is $\Phi_{\rho W}$, where $W$ is the Rademacher variable: $P(W = -1) = P(W = 1) = 1/2$. Let $X = WY$. We are going to show that

$$\mathrm{E}^Z(Y) = 0 \tag{18}$$

$$\mathrm{E}^{X,Z}(Y) = 0 \tag{19}$$

$$\mathrm{E}^Z(X) = \rho Z \tag{20}$$

$$\mathrm{E}^{Y,Z}(X) = Yf(Y,Z), \tag{21}$$

where

$$f(Y,Z) = \frac{\phi_\rho - \phi_{-\rho}}{\phi_\rho + \phi_{-\rho}}(Y,Z); \tag{22}$$

clearly, (10) holds for all $\rho$, but (11) holds if and only if $\rho = 0$.

It follows from the definition of $(Y, Z)$ and $W$:

$$\text{the unconditional joint distribution of } (Y,Z) \text{ is } \frac{1}{2}(\Phi_\rho + \Phi_{-\rho}) \tag{23}$$

$$\text{the marginal distribution of } Y, \text{ as well as } Z, \text{ is } \mathcal{N}(0,1). \tag{24}$$

Let $\phi_\nu^{(1|2)}$ denote the conditional density of $U$ given $V$, where $(U, V) \sim \Phi_\nu$, i.e.,

$$\phi_\nu^{(1|2)}(u|v)\psi(v) = \phi_\nu(u,v). \tag{25}$$

Using the formula $\mathrm{E}^V(U) = \nu V$,

$$\int_\Re u\phi_\nu(u,v)du = \psi(v)\int_\Re u\,\phi_\nu^{(1|2)}(u|v)du = \nu v\psi(v); \tag{26}$$

since $1 - \nu^2 = \mathrm{Var}^V(U)$, using (25) and (14) we obtain

$$\int_\Re u^2\phi_\nu(u,v)du = \psi(v)\bigl(1 - \nu^2 + \nu^2 v^2\bigr). \tag{27}$$

The conditional density of $Y$ given $Z = z$ is obtained by dividing the joint density of $(Y, Z)$ by the marginal density of $Z$; therefore, by (26), (23), and (24), for any $z \in \Re$, the conditional expectation of $Y$ given $Z = z$ is 0, thereby verifying (18).

Since, by the pull-out property, $\mathrm{E}^{X,Z}(Y) = X\mathrm{E}^{X,Z}(W^{-1})$, and $\mathrm{E}(W^{-1}) = 0$, (19) will follow once we show that $(X, Z)$ is independent of $W$. Since $(-Y, Z)$ is a linear transform of $(Y, Z)$, its distribution given $W$ is $\Phi_{-\rho W}$. Consequently, for $B \in \mathcal{B}(\Re^2)$, since $[W = -1] = 1 - [W = 1]$, $P((X,Z) \in B) = \Phi_\rho(B)$ by the averaging property. Therefore, $\mathrm{E}\bigl(\mathrm{E}^W([(X,Z) \in B])A\bigr) = \mathrm{E}(\Phi_\rho(B)A) = \mathrm{E}([(X,Z) \in B]A)$ for $A = \mathfrak{K}$ or $\Omega$, where $\Omega$ is the sample space underlying $W$, $Y$, and $Z$. Since, by the averaging property and the definitions of $W$ and $X$,

$$\mathrm{E}\bigl(\mathrm{E}^W([(X,Z) \in B])[W=1]\bigr) = \mathrm{E}([(X,Z) \in B][W=1]) = \mathrm{E}(\Phi_\rho(B)[W=1])$$

and $\sigma(W) = \{\mathfrak{K}, [W = -1], [W = 1], \Omega\}$, we have $\mathrm{E}^W([(X,Z) \in B]) = \Phi_\rho(B)$. Since

the conditional distribution of $(X, Z)$ given $W$ and the unconditional distribution of $(X, Z)$ both equal $\Phi_\rho$, the independence of $(X, Z)$ and $W$ follows.

Note that (21) will follow once we establish that $\mathrm{E}^{Y,Z}(W) = f(Y, Z)$; since $f$ is a measurable function, it suffices to show, for every $D \in \mathcal{B}(\Re^2)$,

$$\mathrm{E}(W[(Y, Z) \in D]) = \mathrm{E}(f(Y, Z)[(Y, Z) \in D]). \tag{28}$$

Clearly, by the definitions of $W$ and $(Y, Z)$, LHS(28) equals $(\Phi_\rho - \Phi_{-\rho})(D)/2$, which, by (22) and (23), equals RHS(28), establishing (21).

By the chain rule and (21), (20) is equivalent to $\mathrm{E}^Z(Yf(Y, Z)) = \rho Z$, which follows once we show that

$$\mathrm{E}[Yf(Y, Z)[Z \in B]] = \mathrm{E}(\rho Z[Z \in B]) \tag{29}$$

for every $B \in \mathcal{B}(\Re)$. By (22) and (23) again, LHS(29) equals

$$\frac{1}{2} \int\!\!\int_{\Re^2} y[z \in B](\phi_\rho - \phi_{-\rho})(y, z) dy dz,$$

which, using the Fubini Theorem on $\Re^2$, (26) and (24), in that order, equals RHS(29). //

**Remark 4.** For $X \in \mathcal{L}_1$, Remark 2 and Example 1 furnish a counterexample to the conjecture that (10) implies conditional independence of $X$ and $Y$ given $Z$. It turns out that even (10) and (11) combined does not necessarily imply conditional independence of $X$ and $Y$ given $Z$. Consider Example 1 with $\rho = 0$. Clearly, (10) holds, (11) holds by (21), (22), and (20), the conditional distribution of $Y$ given $Z$ is $\mathcal{N}(0, 1)$ by (23), and, since the distribution of $(X, Z)$ is $\Phi_0$ (see the proof of (19) in Example 1), the conditional distribution of $X$ given $Z$ is $\mathcal{N}(0, 1)$ as well. Therefore, conditional independence of $X$ and $Y$ given $Z$ would imply that the distribution of $X+Y$ given $Z$ is $\mathcal{N}(0, 2)$, further implying that the (unconditional) distribution of $X+Y$ is $\mathcal{N}(0, 2)$, which contradicts $2P(X + Y = 0) = 1$ [Majumdar (2017, Example 1)]. //

**Remark 5.** Recall from Example 1 that if $\rho \neq 0$, then (11) does not hold, even though $X \in \mathcal{L}_1$ and both sides of (11) are well-defined. We are going to show that (15) holds in this context, thereby establishing that (15) does not necessarily imply (11) and, by Remark 2, conditional independence of $X$ and $Y$ given $Z$. By (12) and (18) it suffices to show that $\mathrm{E}^Z(YX) = 0$. By the chain rule and (21), $\mathrm{E}^Z(YX) = \mathrm{E}^Z(Y^2 f(Y, Z))$, and it suffices to show $\mathrm{E}(Y^2 f(Y, Z)[Z \in B]) = 0$ for any $B \in \mathcal{B}(\Re)$. Applying (23) and (22),

$$2\mathrm{E}(Y^2 f(Y, Z)[Z \in B]) = \int\!\!\int_{\Re^2} y^2[z \in B](\phi_\rho - \phi_{-\rho})(y, z) dy dz,$$

which, by the Fubini Theorem on $\Re^2$ and (27) applied in that order, equals 0. //

We now present the statement of Theorem 1.

**Theorem 1.** For a pair of conjugate exponents $1 \leq p, q \leq \infty$, let $X \in \mathcal{L}_p$ and $Y \in \mathcal{L}_q$. Then, the following seven conclusions hold:

(i) $\max(p, q) \geq 2$ and $\min(p, q) \leq 2$.

(ii) Conditional mean independence, as in (10), implies zero conditional covariance and
$$\mathrm{E}^{X,Z}(Y) = \alpha(Z) + \beta(Z)X \tag{30}$$
for functions $\alpha, \beta : (\mathfrak{S}, \mathfrak{T}) \mapsto (\Re, \mathcal{B}(\Re))$ such that $\alpha(Z) \in \mathcal{L}_1$ and $\alpha(Z)X \in \mathcal{L}_1$.

(iii) If $p \geq 2$ and (30) holds, then
$$\alpha(Z) = \mathrm{E}^Z(Y) - \beta(Z)\mathrm{E}^Z(X) \tag{31}$$

$$\beta(Z) = \frac{\mathrm{Cov}^Z(X,Y)}{\mathrm{Var}^Z(X)} \left[\mathrm{Var}^Z(X) > 0\right]. \tag{32}$$

(iv) If $p \geq 2$ and (30) holds, then (15) implies (10).

(v) Conditional mean independence, as in (11), implies zero conditional covariance and
$$\mathrm{E}^{Y,Z}(X) = \gamma(Z) + \delta(Z)Y \tag{33}$$
for functions $\gamma, \delta : (\mathfrak{S}, \mathfrak{T}) \mapsto (\Re, \mathcal{B}(\Re))$ such that $\gamma(Z) \in \mathcal{L}_1$ and $\gamma(Z)Y \in \mathcal{L}_1$.

(vi) If $q \geq 2$ and (33) holds, then
$$\gamma(Z) = \mathrm{E}^Z(X) - \delta(Z)\mathrm{E}^Z(Y) \tag{34}$$

$$\delta(Z) = \frac{\mathrm{Cov}^Z(X,Y)}{\mathrm{Var}^Z(Y)} \left[\mathrm{Var}^Z(Y) > 0\right]. \tag{35}$$

(vii) If $q \geq 2$ and (33) holds, then (15) implies (11).

Before presenting the proof of Theorem 1, we show how (30) holds when $X$ is binary, thereby establishing (2), via (31) and (32), for Bernoulli $X, Y \in \mathcal{L}_1$, and arbitrary $Z$.

**Remark 6.** Let $X$ be a binary random variable taking values $b$ and $d$, so that $\mathfrak{X} = \{b, d\}$. Since $X \in \mathcal{L}_\infty$, Theorem 1 will apply for $Y \in \mathcal{L}_1$. For any real valued measurable function $c$ on $\mathfrak{X} \times \mathfrak{S}$, $c(X, Z)$ can be written as
$$c(X, Z) = \frac{dc(b, Z) - bc(d, Z)}{d - b} + \frac{(c(d, Z) - c(b, Z))}{d - b} X.$$

Since $\mathrm{E}^{X,Z}(Y)$ is a real valued $\sigma(X, Z)$ measurable function, there exists a real valued measurable function $e$ on $\mathfrak{X} \times \mathfrak{S}$ such that $\mathrm{E}^{X,Z}(Y) = e(X, Z)$ [Dudley (1989, Theorem 4.2.8)]. Thus, we can write $\mathrm{E}^{X,Z}(Y) = \alpha(Z) + \beta(Z)X$, where $\alpha(Z)$ equals

$$\frac{de(b, Z) - be(d, Z)}{d - b}$$
$$= \frac{(-1)^{\frac{X-d}{b-d}}}{d - b} \left[ Xe\left(b\frac{X - b}{d - b} + d\frac{d - X}{d - b}, Z\right) - (d + b - X)e\left(d\frac{X - b}{d - b} + b\frac{d - X}{d - b}, Z\right) \right].$$

By the triangle inequality,

$$|\alpha(Z)| \leq \frac{\max(|b|, |d|)}{|d - b|} \left[ \left|e\left(b\frac{X - b}{d - b} + d\frac{d - X}{d - b}, Z\right)\right| + \left|e\left(d\frac{X - b}{d - b} + b\frac{d - X}{d - b}, Z\right)\right| \right].$$

By the conditional Jensen's inequality, and the averaging property,

$$\mathrm{E}(|\alpha(Z)|) \leq 2\max(|b|, |d|)(|d - b|)^{-1}\mathrm{E}(|Y|);$$

since $X \in \mathcal{L}_\infty$, (30) holds. By the assertion in (iv) of Theorem 1, (15) implies (10) for Bernoulli $X$, without requiring the propensity score to be non-degenerate. //

Proof of Theorem 1. Since $1 \leq p, q \leq \infty$ and $p^{-1} + q^{-1} = 1$, the assertion in (i) follows.

We only prove the assertions in (ii), (iii), and (iv). Identical arguments, with $X$ and $Y$ interchanged, prove the assertions in (v), (vi), and (vii).

By (13), the chain rule, and the pull-out property,

$$\mathrm{E}^Z(XY) = \mathrm{E}^Z\left(\mathrm{E}^{X,Z}(XY)\right) = \mathrm{E}^Z\left(X\mathrm{E}^{X,Z}(Y)\right). \tag{36}$$

Recall from the Moment inequality that $\mathrm{E}^Z(Y) \in \mathcal{L}_q$. Thus, by the Holder inequality,

$$X\mathrm{E}^Z(Y) \in \mathcal{L}_1. \tag{37}$$

Clearly, (10) implies $\mathrm{RHS}(36) = \mathrm{E}^Z\left(X\mathrm{E}^Z(Y)\right)$, which, by (37) and the pull-out property, equals $\mathrm{E}^Z(X)\mathrm{E}^Z(Y)$, and (15) follows from (12). Further, (30) holds with $\beta(Z) = 0$ and $\alpha(Z) = \mathrm{E}^Z(Y)$; note that, since $q \geq 1$, $\alpha(Z) \in \mathcal{L}_1$ and $\alpha(Z)X \in \mathcal{L}_1$ by (37). That proves the assertion in (ii).

Conversely, assume $p \geq 2$ and (30) holds. If $X$ is $\sigma(Z)$ measurable, then $\mathrm{E}^{X,Z}(Y) = \mathrm{E}^Z(Y)$ and there exists a real valued measurable function $\eta$ on $(\mathfrak{S}, \mathfrak{T})$ such that $X = \eta(Z)$, reducing (30) to $\mathrm{E}^{X,Z}(Y) = \mathrm{E}^Z(Y) = \alpha(Z) + \beta(Z)\eta(Z) = \alpha'(Z)$. In this reparametrization, the "slope" $\beta'$ is equal to 0, and (31) holds for $\alpha'$ by (37). Also, since $X = \mathrm{E}^Z(X)$, we have $\left[\mathrm{Var}^Z(X) > 0\right] = 0$, rendering $\mathrm{RHS}(32) = 0$, thereby vacuously establishing (32) for $\beta'$.

Thus, without loss of generality, we assume that $X$ is not $\sigma(Z)$ measurable, that is, $\mathrm{Var}^Z(X)$ is not almost surely 0. Since $\mathrm{E}^{X,Z}(Y) \in \mathcal{L}_1$ and $\alpha(Z) \in \mathcal{L}_1$, we obtain

$$\beta(Z)X \in \mathcal{L}_1. \tag{38}$$

By the chain rule and the pull-out property, it follows from (30) and (38) that

$$E^Z(Y) = \alpha(Z) + E^Z(\beta(Z)X) = \alpha(Z) + \beta(Z)E^Z(X), \tag{39}$$

whence (31) follows. To show (32), note that (39) implies

$$E^Z(X)E^Z(Y) = \alpha(Z)E^Z(X) + \beta(Z)\left(E^Z(X)\right)^2. \tag{40}$$

By (13) and the pull-out property, $XE^{X,Z}(Y) = E^{X,Z}(XY) \in \mathcal{L}_1$. We can conclude from (30) that $X(\alpha(Z) + \beta(Z)X) \in \mathcal{L}_1$. Since $\alpha(Z)X \in \mathcal{L}_1$, we further obtain that

$$\beta(Z)X^2 \in \mathcal{L}_1. \tag{41}$$

Using (36) with (30), (41), and the pull-out property,

$$E^Z(XY) = E^Z(X(\alpha(Z) + \beta(Z)X)) = \alpha(Z)E^Z(X) + \beta(Z)E^Z(X^2). \tag{42}$$

Subtracting (40) from (42), and using (12) and (14),

$$\text{Cov}^Z(X,Y) = \beta(Z)\text{Var}^Z(X).$$

Multiplying both sides by $k(Z) = [\text{Var}^Z(X) > 0](\text{Var}^Z(X))^{-1}$, we obtain

$$\frac{\text{Cov}^Z(X,Y)}{\text{Var}^Z(X)}[\text{Var}^Z(X) > 0] = \beta(Z)[\text{Var}^Z(X) > 0]. \tag{43}$$

Clearly, (32) will follow from (43) once we show that

$$\beta(Z)[\text{Var}^Z(X) = 0] = 0. \tag{44}$$

By (14) and the averaging property,

$$0 = E(\text{Var}^Z(X)[\text{Var}^Z(X) = 0]) = E\left(\left(X - E^Z(X)\right)^2 [\text{Var}^Z(X) = 0]\right),$$

implying $(X - E^Z(X))^2 [\text{Var}^Z(X) = 0] = 0$. Since the indicator function is idempotent, we obtain $X[\text{Var}^Z(X) = 0] = E^Z(X)[\text{Var}^Z(X) = 0]$; by (39),

$$\alpha(Z)[\text{Var}^Z(X) = 0] + \beta(Z)X[\text{Var}^Z(X) = 0] = E^Z(Y)[\text{Var}^Z(X) = 0],$$

implying $\beta(Z)X[\text{Var}^Z(X) = 0]$ is $\sigma(Z)$ measurable. Since $X$ is not $\sigma(Z)$ measurable, (44) follows, completing the proof of the assertion in (iii).

Finally, assume $p \geq 2$ and (30) holds. By the assertion in (iii), (31) and (32) hold. Then (15) implies $\beta(Z) = 0$, which, in turn, implies, by (30) and (31), that $E^{X,Z}(Y) = E^Z(Y)$, that is, (10) holds. That proves the assertion in (iv). □

**Remark 7.** To appreciate the power of Theorem 1, let us review how a direct proof of (2) for Bernoulli $X$, that makes no use of Theorem 1 and Remark 6, might work. As discussed in Section 1, the first order of business here is to show that $g(X,Z) \in \mathcal{L}_1$. We can write, using (12), (14), and the idempotence of $X$,

$$g(X,Z) = \frac{\mathrm{E}^Z\big((X - \mathrm{E}^Z(X))Y\big)}{\sqrt{\mathrm{E}^Z(X)(1 - \mathrm{E}^Z(X))}}\big[0 < \mathrm{E}^Z(X) < 1\big] \frac{(X - \mathrm{E}^Z(X))}{\sqrt{\mathrm{E}^Z(X)(1 - \mathrm{E}^Z(X))}}. \quad (45)$$

It is reasonably clear from the three factor representation of $g(X,Z)$ in (45) that $g(X,Z) \notin \mathcal{L}_\infty$; therefore, we have work to do to show that $g(X,Z) \in \mathcal{L}_1$.

Since the indicator function is idempotent, writing

$$\frac{(X - \mathrm{E}^Z(X))}{\sqrt{\mathrm{E}^Z(X)(1 - \mathrm{E}^Z(X))}}\big[0 < \mathrm{E}^Z(X) < 1\big] = h(X,Z) \quad (46)$$

and using the pull-out property, we can obtain

$$g(X,Z) = h(X,Z)\mathrm{E}^Z(h(X,Z)Y), \quad (47)$$

once we show that

$$h(X,Z)Y \in \mathcal{L}_1. \quad (48)$$

Since $h(X,Z) \notin \mathcal{L}_\infty$, (48) cannot be obtained from the assumption $Y \in \mathcal{L}_1$.

We now show that

$$\mathrm{E}\big(h^2(X,Z)\big) \leq 1. \quad (49)$$

By Theorem 2.10 of Folland (1999), there exists an increasing sequence $\{W_n : n \geq 1\}$ of nonnegative simple functions on $(\Omega, \sigma(Z))$ such that

$$W_n \text{ converges pointwise to } \frac{[0 < \mathrm{E}^Z(X) < 1]}{\mathrm{E}^Z(X)(1 - \mathrm{E}^Z(X))}. \quad (50)$$

By (46) and (50), $(X - \mathrm{E}^Z(X))^2 W_n$ is an increasing sequence of nonnegative functions converging pointwise to $h^2(X,Z)$; by the Monotone Convergence Theorem (MCT),

$$\mathrm{E}\big(h^2(X,Z)\big) = \lim_{n\to\infty} \mathrm{E}\Big((X - \mathrm{E}^Z(X))^2 W_n\Big). \quad (51)$$

For every $n \geq 1$, since $W_n$ is a simple function on $(\Omega, \sigma(Z))$ and $|X - \mathrm{E}^Z(X)| \leq 1$, $(X - \mathrm{E}^Z(X))^2 W_n$ is bounded; applying the pull-out property,

$$\mathrm{E}^Z\Big((X - \mathrm{E}^Z(X))^2 W_n\Big) = W_n \mathrm{Var}^Z(X). \quad (52)$$

By (51), (52), and the averaging property, $\mathrm{E}(h^2(X,Z)) = \lim_{n\to\infty} \mathrm{E}\big(W_n \mathrm{Var}^Z(X)\big)$; by (50), $W_n \mathrm{Var}^Z(X)$ is an increasing sequence of nonnegative functions converging pointwise to $[0 < \mathrm{E}^Z(X) < 1]$, whence (49) follows by the MCT.

From (47), by the Cauchy-Schwartz (C-S hereinafter) inequality and (49),

$$\left(\mathrm{E}(|g(X,Z)|)\right)^2 \leq \mathrm{E}\left(\left(\mathrm{E}^Z(h(X,Z)Y)\right)^2\right). \tag{53}$$

We can think of two ways of concluding $g(X,Z) \in \mathcal{L}_1$ from (53).

If $Y \in \mathcal{L}_\infty$, (49) implies $h(X,Z)Y \in \mathcal{L}_2$, whence (48), and hence (47), follow by the Moment inequality; further, by (53), the conditional Jensen's inequality, and the averaging property, $\|g(X,Z)\|_1 \leq \|h(X,Z)Y\|_2$.

If $Y \in \mathcal{L}_2$, (48), and hence (47), follow from (49) by the C-S inequality. By (49), (46) and the pull-out property, $\mathrm{E}^Z(h^2(X,Z)) = [0 < \mathrm{E}^Z(X) < 1] \leq 1$. "Since"

$$\left(\mathrm{E}^Z(h(X,Z)Y)\right)^2 \leq \mathrm{E}^Z\left(h^2(X,Z)\right)\mathrm{E}^Z\left(Y^2\right) \leq \mathrm{E}^Z\left(Y^2\right), \tag{54}$$

(53) implies $\|g(X,Z)\|_1 \leq \|Y\|_2$. However, the first inequality in (54) "follows" by applying the C-S inequality in $\mathcal{L}_2\left(\{0,1\} \times \Re \times \mathfrak{S}, 2^{\{0,1\}} \otimes \mathcal{B}(\Re) \otimes \mathfrak{T}, P_{(X,Y,Z)|Z}\right)$ to the functions $f_1$ and $f_2$, where $f_1(x,y,z) = h(x,z)$ and $f_2(x,y,z) = y$. For that argument to be valid, $P_{(X,Y,Z)|Z}$ has to exist. However, the existence of $P_{(X,Y,Z)|Z}$ is not guaranteed for arbitrary $\mathfrak{S}$ (assuming $\mathfrak{S}$ to be a Polish space would suffice).

Even after we make one of the additional assumptions and establish $g(X,Z) \in \mathcal{L}_1$, non-trivial work remains to be done to verify (2), reiterating that the right path to (2) for Bernoulli $X$ is through Theorem 1 and Remark 6. //

**3. Conclusion.** Using the well-understood properties of the conditional distribution of one component of a multivariate Normal vector given another component, we showed that if $(X,Y,Z)$ is multivariate Normal, where $X, Y \in \Re$, then the formula in (2) holds. We think that is a better explanation of why zero conditional covariance implies conditional mean independence in the multivariate Normal case, compared to the sledgehammer explanation that zero conditional covariance implies conditional independence in that case and conditional independence implies conditional mean independence in general. Further, that leads to the conjecture that the real reason why zero conditional covariance implies conditional mean independence when $X$ is Bernoulli is because the formula in (2) holds in that case as well. We settled that conjecture via Theorem 1 and Remark 6. We argued in Remark 7 that a direct verification of the formula in (2) for Bernoulli $X$ is not possible without additional moment assumptions on $Y$: if one is willing to assume $Y \in \mathcal{L}_\infty$, then one can leave $Z$ arbitrary, but if one wants to only assume $Y \in \mathcal{L}_2$, then further structural assumptions on $Z$ are needed.

**Acknowledgment.** I would like to thank Professor Jungbin Hwang, whose graduate econometrics seminar was instrumental in getting me started on this research, and Professor Jeffrey Ladewig, who introduced me to the seminal work of Anthony Downs. Professor Hwang carefully read multiple drafts to help me improve the exposition.

RAJESHWARI MAJUMDAR
rajeshwari.majumdar@uconn.edu
PO Box 47
Coventry, CT 06238